\documentclass[12pt,twoside,a4paper,reqno]{amsart}

\usepackage{graphicx}
\usepackage{graphics}
\usepackage{amssymb}
\usepackage{eucal}
\usepackage{amsmath}
\usepackage{amsthm}
\newtheorem{theorem}{Theorem}[section]
\newtheorem{lemma}[theorem]{Lemma}

\newtheorem{corollary}[theorem]{Corollary}

\theoremstyle{definition}
\newtheorem{definition}[theorem]{Definition}

\theoremstyle{remark}

\numberwithin{equation}{section}

\begin{document}

\title{Some points on $c$-$K$-$g$-frames and their duals}

\author{Morteza Rahmani}
\address{Young Researchers and Elite Club, Ilkhchi Branch, Islamic Azad University,
Ilkhchi, Iran
}
\email{morteza.rahmany@gmail.com}

\author{Esmaeil Alizadeh }
\address{Department of Mathematics, Marand Branch, Islamic Azad University, Marand, Iran}
\email{e\_alizadeh@marandiau.ac.ir}

\subjclass[2010]{Primary 42C15, 46C05}



\keywords{$cg$-frame, $c$-$K$-$g$-frame, Dual   $c$-$K$-$g$-frame, Atomic $cg$-system}

\begin{abstract}
In this paper we study some new properties of $c$-$K$-$g$-frames in a Hilbert space $H$. We study  duals of $c$-$K$-$g$-frames
and give some characterizations of $c$-$K$-$g$-frames and their duals. Also, we verify the relationships between $c$-$K$-$g$-frames  and atomic $cg$-systems.
Precisely, we show that these two concepts are equivalent. Finally, we find some new atomic $cg$-systems from given ones.
\end{abstract}

\maketitle

\section{Introduction}

Frames  in Hilbert spaces were introduced by Duffin and Schaeffer \cite{Duff} in 1952 to study some deep problems in nonharmonic Fourier series.
A frame for a Hilbert space $H$ provides a linear combination of the elements of frame for each element in $H$,
 but linear independence between
the frame elements is not required. In other word, a frame can be thought
as a basis to which one has added more elements.
Motivated by the theory of coherent states, this concept was generalized  to families indexed by some locally compact space endowed with a Radon measure. This approach leads to the notion of continuous frames. More details of these kinds of frames is discussed in \cite{1, 2, 14, Rah}.

Gavruta  introduced the $K$-frames in Hilbert spaces to study atomic
decomposition systems  and discussed some properties of them \cite{k-frame}.
Afterward, $K$-$g$-frames have been introduced in \cite{10}. Also, the concept of continuous
$K$-$g$-frames is introduced in \cite{ckg}.

Throughout this paper, $(\Omega,\mu)$ is a measure space with positive measure $\mu$,  $H$, $H_1$, $H_2$ and $H_{\omega}$ are separable Hilbert spaces and $B(H,H_{\omega})$ is the set of all bounded linear operators from $H$ into $H_{\omega}$, $\omega\in\Omega$. If $H_{\omega}=H$, then $B(H,H)$ will be denoted by $B(H)$.

\begin{lemma}{\rm(\cite{dag})} \label{1.1}
Let $L_1\in B(H_1,H)$ and $L_2\in B(H_2,H)$. Then the following assertions are equivalent:
\begin{enumerate}
\item[(i)]
$R(L_1)\subseteq R(L_2)$.
\item[(ii)]
$L_1L_1^*\leq\lambda^2 L_2L_2^*$ for some $\lambda>0$.
\item[(iii)]
There exists an operator  $U\in B(H_1,H_2)$ such that $L_1=L_2 U$.
\end{enumerate}
Moreover, if $(i)$, $(ii)$ and $(iii)$ are valid, then there exists a unique operator
$U$ so that
\begin{enumerate}
\item[(i)]
$\|U\|^2=inf \big\{\mu:  L_1L_1^*\leq \mu L_2L_2^* \big\}$,
\item[(ii)]
$N(L_1)=N(U)$,
\item[(iii)]
$R(U)\subseteq \overline{R(L_2)^*}$.
\end{enumerate}
\end{lemma}

\begin{definition}
Let $\varphi\in\Pi_{\omega\in \Omega}H_\omega$. We say that $\varphi$ is  strongly measurable  if $\varphi$ as
a mapping of $\Omega$ to $\oplus_{\omega \in \Omega}H_{\omega}$ is measurable, where
$$\Pi_{\omega\in \Omega}H_\omega=\big\{f:\Omega\longrightarrow \cup_{\omega\in \Omega} H_\omega~;~f(\omega) \in H_\omega\big\}.$$
\end{definition}

\begin{definition}
Consider the set
\begin{align*}
\big(\oplus_{\omega \in \Omega}H_{\omega},\mu\big)_{L^2}=
\big\{F |~ F ~is~strongly~measurable, \int_{\Omega}\|F(\omega)\|^{2}d\mu(\omega) < \infty\big\},
\end{align*}
with inner product given by
\[ \langle F, G \rangle=\int_{\Omega} \langle F(\omega), G(\omega) \rangle d\mu(\omega).\]
It can be proved that $\Big(\oplus_{\omega\in\Omega}H_{\omega}, \mu\Big)_{L^{2}}$ is a Hilbert space (\cite{1}). We will denote the norm of $F\in \Big(\oplus_{\omega\in\Omega}H_{\omega}, \mu\Big)_{L^{2}}$ by $\|F\|_{2}$.
\end{definition}

We review the definition of continuous $g$-frames.

\begin{definition}
We call $\{\Lambda_{\omega}\in B(H,H_\omega): \omega\in \Omega \}$ a  continuous generalized frame, or simply a   $cg$-frame, for $H$ with respect to $\{H_\omega\}_{\omega\in \Omega}$  if:
\begin{enumerate}
\item[(i)]
 for each $f\in H$, $\{\Lambda_{\omega}f\}_{\omega\in \Omega}$ is strongly measurable,
\item[(ii)]
 there are two positive constants $A$ and $B$ such that
\begin{align}\label{ccg}
A\|f\|^{2}\leq\int_{\Omega}\|\Lambda_\omega f\|^{2}d\mu(\omega)\leq B\|f\|^{2},\quad f\in H.
\end{align}
\end{enumerate}
We call $A$ and $B$ the lower and upper  $cg$-frame bounds, respectively. If $A,B$ can be chosen such that $A=B$, then $\{\Lambda_{\omega}\}_{\omega\in \Omega}$ is called a  tight $cg$-frame  and if
$A = B = 1$, it is called a Parseval $cg$-frame. A family $\{\Lambda_{\omega}\}_{\omega\in \Omega}$ is called  $cg$-Bessel family  if the second inequality in (\ref{ccg}) holds.
\end{definition}

\begin{theorem}{\rm(\cite{1})}
Let $\{\Lambda_{\omega}\}_{\omega\in\Omega}$ be a  $cg$-Bessel family for $H$ with respect to $\{H_{\omega}\}_{\omega\in \Omega}$ with   bound $B$. Then the mapping  $T$ of  $\Big(\oplus_{\omega\in\Omega}H_{\omega}, \mu\Big)_{L^{2}}$ to $H$ weakly defined by
\[\langle TF, g \rangle=\int_{\Omega}\langle \Lambda_{\omega}^{*}F(\omega), g \rangle d\mu(\omega), \quad F\in \Big(\oplus_{\omega\in\Omega}H_{\omega}, \mu\Big)_{L^{2}},~g\in H,\]
is linear and bounded with $\|T\|\leq\sqrt{B}$. Furthermore for each $g\in H$ and $\omega\in \Omega$,
\[T^{*}(g)(\omega)=\Lambda_{\omega}g.\]
\end{theorem}
The operator $T$ is called the synthesis operator of $\{\Lambda_{\omega}\}_{\omega\in\Omega}$ and $T^*$ is called the analysis operator of $\{\Lambda_{\omega}\}_{\omega\in\Omega}$.

\vspace{1.5mm}
The continuous version of $K$-$g$-frames have been introduced in \cite{ckg} as following:
\begin{definition}
Let $K\in B(H)$. A family $\Lambda=\{\Lambda_{\omega}\in B(H, H_{\omega}):~\omega\in\Omega\}$ is called a continuous $K$-$g$-frame, or $c$-$K$-g-frame, for $H$ with respect to $\{H_{\omega}\}_{\omega\in \Omega}$  if:
\begin{enumerate}
\item[(i)]
$\{\Lambda_{\omega}f\}_{\omega\in\Omega}$ is strongly measurable for each $f\in H$,
\item[(ii)]
there exist constants  $0<A\leq B<\infty$ such that
\begin{equation}\label{e}
A\|K^{*}f\|^{2}\leq \int_{\Omega} \|\Lambda_{\omega}f\|^{2}\,d\mu(\omega)\leq B\|f\|^{2},\quad f\in H.
\end{equation}
\end{enumerate}
\end{definition}
The constants $A$, $B$ are called lower and upper $c$-$K$-$g$-frame bounds, respectively.
If $A$, $B$ can be chosen such that $A=B$, then $\{\Lambda_{\omega}\}_{\omega\in\Omega}$ is called a tight $c$-$K$-$g$-frame and if $A=B=1$, it is called a Parseval $c$-$K$-$g$-frame. The family $\{\Lambda_{\omega}\}_{\omega\in\Omega}$ is called a $c$-$g$-Bessel family if the right hand inequality in (\ref{e}) holds. In this case, $B$ is called the Bessel constant.

\vspace{1.5mm}
Now, suppose that $\{\Lambda_{\omega}\}_{\omega\in\Omega}$ is a $c$-$K$-$g$-frame for $H$ with respect to $\{H_{\omega}\}_{\omega\in\Omega}$ with frame bounds $A$, $B$. The $c$-$K$-$g$-frame operator $S:H\longrightarrow H$ is weakly defined by
\begin{align*}
\langle Sf, g\rangle= \int_{\Omega} \langle f, \Lambda_{\omega}^{*} \Lambda_{\omega}g \rangle\, d\mu(\omega),\quad  f, g\in H.
\end{align*}
Therefore,
$$AKK^*\leq S\leq BI.$$
\begin{lemma}{\rm(\cite{ckg})} \label{x}
Let $\{\Lambda_{\omega}\}_{\omega\in \Omega}$ be a $cg$-Bessel family  for $H$ with respect to $\{H_{\omega}\}_{\omega\in\Omega}$. Then $\{\Lambda_{\omega}\}_{\omega\in\Omega}$ is a $c$-$K$-$g$-frame for $H$ with respect to $\{H_{\omega}\}_{\omega\in\Omega}$ if and only if there exists a  constant $A>0$ such that $S\geq AKK^{*}$, where $S$ is the frame operator of $\{\Lambda_{\omega}\}_{\omega\in \Omega}$.
\end{lemma}

Duals of $c$-$K$-$g$-frames have been introduced in \cite{sckg} as following:
\begin{definition}
Let $\Lambda=\{\Lambda_{\omega}\}_{\omega\in\Omega}$ be a $c$-$K$-$g$-frame for $H$ with respect to $\{H_{\omega}\}_{\omega\in\Omega}$. A $cg$-Bessel family $\Gamma=\{\Gamma_{\omega}\}_{\omega\in\Omega}$ for $H$ is called a dual $c$-$K$-$g$-Bessel family of $\Lambda$ if for each $f, h\in H$,
$$\langle Kf, h\rangle=\int_{\Omega}\langle\Lambda^*_{\omega}\Gamma_{\omega}f, h\rangle\,d\mu(\omega).$$
\end{definition}
In this case,  $\Lambda$ and $\Gamma$ are are called pair duals.
\section{Characterizing dual $c$-$K$-$g$-frames}
In this is section,  we study the duals of $c$-$K$-$g$-frame and verify their characterizations.

The following theorem characterizes a $c$-$K$-$g$-frame by operator decompositions and also gives a sufficient condition for a $cg$-Bessel family to be dual   of a $c$-$K$-$g$-Bessel family and a $c$-$K$-$g$-frame.
\begin{theorem}\label{a1}
Let  $K\in B(H)$ and  $\{\Lambda_{\omega}\}_{\omega\in\Omega}$ be a $cg$-Bessel family for $H$ with synthesis operator $T_{\Lambda}$. Then  $\{\Lambda_{\omega}\}_{\omega\in\Omega}$ is a c-$K$-g-frame for $H$ if and only if there exists a bounded operator $\Phi: \Big(\oplus_{\omega\in \Omega}H_{\omega}, \mu\Big)_{L^{2}}\longrightarrow H$ such that  $K^*=\Phi T^*_{\Lambda}$. Furthermore, if  $\Gamma=\{\Gamma_{\omega}\}_{\omega\in\Omega}$ where
$$\Gamma_{\omega}g=(\Phi^*g)(\omega), \quad g\in H,~ \omega\in\Omega,$$
then  $\{\Gamma_{\omega}\}_{\omega\in\Omega}$ is a dual $c$-$K$-$g$-Bessel family of  $\{\Lambda_{\omega}\}_{\omega\in\Omega}$.
\begin{proof}
Let $\{\Lambda_{\omega}\}_{\omega\in\Omega}$ be a $c$-$K$-$g$-frame for $H$. Then by  Theorem 2.5 in \cite{ckg}, $R(K)\subseteq R(T_{\Lambda})$. By Lemma \ref{1.1}, there exists a bounded operator
$$\Phi: \Big(\oplus_{\omega\in \Omega}H_{\omega}, \mu\Big)_{L^{2}}\longrightarrow H$$
such that $K=T_{\Lambda}\Phi^*$. So  $K^*=\Phi T^*_{\Lambda}$.

\vspace{1.5mm}
Now for the opposite implication, assume that  there exists a bounded operator $\Phi: \Big(\oplus_{\omega\in \Omega}H_{\omega}, \mu\Big)_{L^{2}}\longrightarrow H$ such that $K^*=\Phi T^*_{\Lambda}$. Then  $K=T_{\Lambda}\Phi^*$ and by lemma \ref{1.1}, $R(K)\subseteq R(T_{\Lambda})$, also  by  Theorem 2.5 in \cite{ckg}, $\{\Lambda_{\omega}\}_{\omega\in\Omega}$ is a $c$-$K$-$g$-frame for $H$.
For each $f,g\in H$,  we have
\begin{align*}
\langle K^*f, g\rangle&=\langle \Phi T^*_{\Lambda}f, g\rangle=\langle T^*_{\Lambda}f, \Phi^*g\rangle\\
&=\int_{\Omega}\big\langle\Lambda_{\omega}f, (\Phi^*g)(\omega)\rangle\,d\mu(\omega)\\
&=\int_{\Omega}\big\langle\Lambda_{\omega}f,\Gamma_{\omega}g\rangle\,d\mu(\omega)\\
&=\int_{\Omega}\big\langle\Gamma^*_{\omega}\Lambda_{\omega}f,g\rangle\,d\mu(\omega),
\end{align*}
 where
$\Gamma_{\omega}g=(\Phi^*g)(\omega)$, $g\in H,~ \omega\in\Omega$.
Also
$$\int_{\Omega}\Vert\Gamma_{\omega}g\Vert^2 d\,\mu=\int_{\Omega}\Vert(\Phi^*g)(\omega)\Vert^2 d\,\mu=\|\Phi^*g\|^{2}_{2}\leq\|\Phi\|^{2}_{2}\|g\|^{2}, \quad g\in H.$$
So $\Gamma=\{\Gamma_{\omega}\}_{\omega\in\Omega}$ is a dual $c$-$K$-$g$-Bessel family of  $\{\Lambda_{\omega}\}_{\omega\in\Omega}$.
\end{proof}
\end{theorem}
\begin{theorem}\label{a2}
Assume that $K\in B(H)$ and  $\{\Lambda_{\omega}\}_{\omega\in\Omega}$ is a $c$-$K$-$g$-frame for $H$ with synthesis operator $T_{\Lambda}$. Then $\{\Gamma_{\omega}\}_{\omega\in\Omega}$ is a dual $c$-$K$-$g$-Bessel family of  $\{\Lambda_{\omega}\}_{\omega\in\Omega}$ if and only if there exists a bounded operator
$$\Phi: \Big(\oplus_{\omega\in \Omega}H_{\omega}, \mu\Big)_{L^{2}}\longrightarrow H$$
such that $K^*=\Phi T^*_{\Lambda}$ and $\Gamma_{\omega}g=(\Phi^*g)(\omega), ~g\in H,~\omega\in\Omega$.
\end{theorem}
\begin{proof}
Suppose that  $\{\Gamma_{\omega}\}_{\omega\in\Omega}$  is a dual  $c$-$K$-$g$-Bessel family of $\{\Lambda_{\omega}\}_{\omega\in\Omega}$. Consider $\Phi$ as $\Phi=T_{\Gamma}$. Then
$\Phi^*g(\omega)=\Gamma_\omega g$, $g\in H$,  $\omega\in\Omega$.
By assumption, for each $f,g\in H$,
\begin{align*}
\langle K^*f, g\rangle =&\int_{\Omega} \langle \Gamma^*_{\omega}\Lambda_{\omega}f,g\rangle d\mu(\omega)
=\langle T_{\Gamma}\{\Lambda_{\omega}f\}_{\omega\in \Omega}, g\rangle
\\=&\langle \Phi\{\Lambda_{\omega} f\}_{\omega\in\Omega},g\rangle
=\langle \Phi T^*_{\Lambda}f,g\rangle.
\end{align*}
Therefore, $K^*=\Phi T^*_{\Lambda}$.

\vspace{1mm}
The converse  implication has been proved in Theorem \ref{a1}.
\end{proof}

\begin{theorem}
Let  $K\in B(H)$ and $\{\Lambda_{\omega}\}_{\omega\in\Omega}$ be a c-$K$-g-frame for $H$ with optimal lower bound $A$. Suppose that  $\Gamma=\{\Gamma_{\omega}\}_{\omega\in\Omega}$ is a dual c-$K$-g-Bessel family of $\{\Lambda_{\omega}\}_{\omega\in\Omega}$. Then $ \|T_{\Gamma}\|^{2}\geq \frac{1}{A}$, where $T_{\Gamma}$ is the synthesis operator of $\{\Gamma_{\omega}\}_{\omega\in\Omega}$. Furthermore, there exists a unique  dual c-$K$-g-Bessel family $\Theta=\{\Theta_{\omega}\}_{\omega\in \Omega}$   such that $\|T_{\Theta}\|^{ 2}=A$, where $T_{\Theta}$ is the synthesis operate of $\Theta$.
\end{theorem}
\begin{proof}
Let $A$ is the optimal lower $c$-$K$-$g$-frame bound of  $\{\Lambda_{\omega}\}_{\omega\in\Omega}$, then for each $f\in H$,
$$A\|K^*f\|^{2}\leq\int_{\Omega}\Vert\Lambda_{\omega}f\Vert^2 d\,\mu(\omega).$$
Then
$$\|K^*f\|^{2}\leq\frac{1}{A}\|T^*_{\Lambda}f\|^{2},  \quad f\in H.$$
By Theorem 3.1 in \cite{sckg}, we have $K=T_{\Lambda}T^*_{\Gamma}$.
So for each $f\in H$, we obtain
\begin{align*}
\|K^*f\|^{2}&=\langle K^*f, K^*f\rangle=\langle KK^*f, f\rangle=\langle T_{\Lambda}T^*_{\Gamma}T_{\Gamma}T^*_{\Lambda}f, f\rangle\\
&=\langle T_{\Gamma}T^*_{\Lambda}f, T_{\Gamma}T^*_{\Lambda}f\rangle=\|T_{\Gamma}T^*_{\Lambda}f\|^{2}\leq\|T_{\Gamma}\|^{2}\|T^*_{\Lambda}f\|^{2},
\end{align*}
that is,
$ \frac{1}{\|T_\Gamma\|^{2}}\|K^*f\|^{2}\leq\|T^*_{\Lambda}f\|^{2}.$
Since
\begin{align*}
A&=sup\big\{\lambda>0: \lambda\|K^*f\|^{2}\leq\|T^*_{\Lambda}f\|^{2},~ f\in H\big\}\\
&=inf \big\{\mu: \|K^*f\|^{2}\leq\mu\|T^*_{\Lambda}f\|^{2},~ f\in H\big\},
\end{align*}
So $\|T_\Gamma\|^{2}\geq \frac{1}{A}.$
By Theorem 2.5 in \cite{ckg}, $R(K)\subseteq R(T_{\Lambda})$, so by lemma \ref{1.1}, there exists a unique bounded operator $\Phi: \Big(\oplus_{\omega\in \Omega}H_{\omega}, \mu\Big)_{L^{2}}\longrightarrow H$  such that $K^*=\Phi T^*_{\Lambda}$ and
$$ \|\Phi\|^{2}=inf \big\{\mu:  \|K^*f\|^{2}\leq\mu\|T^*_{\Lambda}f\|^{2},~ f\in H\big\}=A$$
Let $\{\Theta_{\omega}\}_{\omega\in\Omega}$ be the family which for each $\omega\in\Omega$, $\Theta_{\omega}$ is defined by $$\Theta_{\omega}f=(\Phi^*f)(\omega), \quad  f\in H.$$ So by Theorem \ref{a2}, $\{\Theta_{\omega}\}_{\omega\in\Omega}$ is a  dual c-$K$-g-frame for $\{\Lambda_{\omega}\}_{\omega\in\Omega}$.
For each $f\in H$  and $\omega\in\Omega$,
$$(T^*_{\Theta}f)(\omega)=\Theta_{\omega}f=(\Phi^*f)(\omega).$$
So $T^*_{\Theta}=\Phi^*$ and $\|T_{\Theta}\|^{2}=A$.
\end{proof}
The following theorem is a generalization of Theorem 2.6 in \cite{Li} for  continuous version.
\begin{theorem}
Assume that $K \in B(H)$ is with closed range and  $\{\Lambda_{\omega}\}_{\omega\in\Omega}$ is a $cg$-Bessel family for $H$ with the frame operator $S_{\Lambda}$. If $\{\Lambda_{\omega}\}_{\omega\in\Omega}$ has a dual $cg$-frame on $R(K)$ and $S_{\Lambda}(R(K))\subseteq R(K)$, then it is a $c$-$K$-$g$-frame for $H$.
\end{theorem}
\begin{proof}
Suppose  that $\{\Gamma_{\omega}\}_{\omega\in\Omega}$ is a dual $cg$-frame of $\{\Lambda_{\omega}\}_{\omega\in\Omega}$ on $R(K)$. For each  $f\in H$, we can write $f=f_{1}+f_{2}$, where $ f_{1}\in R(K)$ and $f_{2}\in (R(K))^\perp$. Thus
\begin{align*}
\int_{\Omega}\Vert\Lambda_{\omega}f\Vert^2 d\,\mu(\omega)&=\int_{\Omega}\Vert\Lambda_{\omega}(f_{1}+f_{2})\Vert^2 d\,\mu(\omega)\\
&=\int_{\Omega}\Vert\Lambda_{\omega}f_{1}\Vert^2 d\,\mu(\omega)+\int_{\Omega}\Vert\Lambda_{\omega}f_{2}\Vert^2 d\,\mu(\omega)\\
&+2Re\int_{\Omega}\big\langle\Lambda^*_{\omega}\Lambda_{\omega}f_{1},f_{2}\rangle\,d\mu(\omega).
\end{align*}
Since $S_{\Lambda}f_{1}\in S_{\Lambda}(R(K))\subseteq R(K)$, we have
\begin{align*}
 \int_{\Omega}\big\langle\Lambda^*_{\omega}\Lambda_{\omega}f_{1},f_{2}\rangle\,d\mu(\omega)=\langle S_{\Lambda}f_{1},f_{2}\rangle=0.
\end{align*}
Hence
$$\int_{\Omega}\Vert\Lambda_{\omega}f\Vert^2 d\,\mu(\omega)=\int_{\Omega}\Vert\Lambda_{\omega}f_{1}\Vert^2 d\,\mu(\omega)+\int_{\Omega}\Vert\Lambda_{\omega}f_{2}\Vert^2 d\,\mu(\omega), \quad f\in H.$$
Note that $ker(K^*)=(R(K))^\perp$ and by the definition of  dual $cg$-frames in \cite{1}, for each $f\in H$, we have
\begin{align*}
\|K^*f\|^{2}&=\|K^*(f_{1}+f_{2})\|^{2}=\|K^*f_{1}\|^{2}=\vert\langle K^*f_{1}, K^*f_{1}\rangle\vert=\vert\langle KK^*f_{1}, f_{1}\rangle\vert\\
&=\Big|\int_{\Omega}\big\langle\Gamma_{\omega}KK^*f_{1},\Lambda_{\omega}f_{1}\rangle\,d\mu(\omega)\Big|\\
&\leq\int_{\Omega}\|\Gamma_{\omega}KK^*f_{1}\|\|\Lambda_{\omega}f_{1}\|\,d\mu(\omega)\\
&\leq\Big(\int_{\Omega}\|\Gamma_{\omega}KK^*f_{1}\|^{2}d\mu(\omega)\Big)^{\frac{1}{2}}\Big(\int_{\Omega}\|\Lambda_{\omega}f_{1}\|^{2}d\mu(\omega)\Big)^{\frac{1}{2}}\\
&\leq(B\|K\|^{2}\|K^{*}f_{1}\|^{2})^{\frac{1}{2}}\Big(\int_{\Omega}\|\Lambda_{\omega}f_{1}\|^{2}d\mu(\omega)\Big)^{\frac{1}{2}},
\end{align*}
where $B$ is the $cg$-Bessel bound of $\{\Gamma_{\omega}\}_{\omega\in\Omega}$, so
$$\frac{1}{B\|K\|^{2}}\|K^*f\|^{2}\leq\int_{\Omega}\|\Lambda_{\omega}f_{1}\|^{2}d\mu(\omega),$$
 and hence
\begin{align*}
\int_{\Omega}\Vert\Lambda_{\omega}f\Vert^2 d\,\mu(\omega)&=\int_{\Omega}\Vert\Lambda_{\omega}f_{1}\Vert^2 d\mu(\omega)+\int_{\Omega}\Vert\Lambda_{\omega}f_{2}\Vert^2 d\mu(\omega)\\
&\geq\int_{\Omega}\Vert\Lambda_{\omega}f_{1}\Vert^2 d\mu(\omega)\\
&\geq\frac{1}{B\|K\|^{2}}\|K^*f\|^{2}.
\end{align*}
\end{proof}
\section{Continuous atomic  $g$-systems and $c$-$K$-$g$-frames}
In this section, we study  the properties of continuous  atomic $g$-systems for an operator $K\in B(H)$ and verify the relationship between this concept and $c$-$K$-$g$-frames.
The continuous version of atomic systems for a family of operator is defined as below (\cite{nob}):
\begin{definition}
Suppose that $K\in B(H)$. A family $\Lambda=\{\Lambda_{\omega}\in B(H, H_{\omega}):~\omega\in\Omega\}$ is called a continuous atomic  $g$-system for $K$, or simply  an  atomic $cg$-system for $K$, if the following conditions hold:
\begin{enumerate}
\item[(i)]
$\{\Lambda_{\omega}\}_{\omega\in\Omega}$ is  a $cg$-Bessel family,
\item[(ii)]
there exists a constant $C>0 $ such that for each $ f\in H$, there exists a $\varphi\in \Big(\oplus_{\omega\in \Omega}H_{\omega}, \mu\Big)_{L^{2}}$   such that $\|\varphi\|_{2}\leq C\|f\|$ and for each $ g\in H$,
$$\langle Kf, g\rangle=\int_{\Omega}\langle\Lambda^*_{\omega}\varphi_{f}(\omega), g\rangle\,d\mu(\omega).$$
\end{enumerate}
\end{definition}
Now, we present a characterization for  atomic $cg$-systems.
\begin{theorem}\label{3.2}
Let  $\{\Lambda_{\omega} \in B(H, H_{\omega}):~\omega\in\Omega\}$  be a family of linear operators. Then the following statements are equivalent.
\begin{enumerate}
\item[(i)]
$\{\Lambda_{\omega}\}_{\omega\in\Omega}$ is  an atomic $cg$-system for $K$.
\item[(ii)]
$\{\Lambda_{\omega}\}_{\omega\in\Omega}$ is $c$-$K$-g-frame for $H$.
\item[(iii)]
There exists a $cg$-Bessel family $\{\Gamma_{\omega}\}_{\omega\in\Omega}$ for $H$ with respect to $\{H_{\omega}\}_{\omega\in\Omega}$ such that
\begin{equation}\label{m}
\langle Kf, h \rangle =\int_{\Omega}\langle\Lambda_{\omega}^{*}\Gamma_{\omega}f, h \rangle d\mu(\omega), \quad f, h\in H.
\end{equation}
\end{enumerate}
\end{theorem}
\begin{proof}
$ (i)\Rightarrow (ii)$  There exists a $C>0$ such that  for each $h\in H$, there exists a $\varphi\in \Big(\oplus_{\omega\in \Omega}H_{\omega}, \mu\Big)_{L^{2}}$   so  that $\|\varphi\|_{2}\leq C\|K^*h\|$. Therefore for each $h\in H$,
\begin{align*}
\|K^*h\|^{2}&=\vert\langle KK^*h, h\rangle\vert
=\Big\vert\int_{\Omega}\langle\Lambda^*_{\omega}\varphi(\omega), h\rangle\,d\mu(\omega)\Big\vert\\
&=\Big\vert\int_{\Omega}\langle\varphi(\omega),\Lambda_{\omega} h\rangle\,d\mu(\omega)\Big\vert
\leq\int_{\Omega}\|\varphi(\omega)\|\|\Lambda_{\omega}h\|d\mu(\omega)\\
&\leq\Big(\int_{\Omega}\|\varphi(\omega)\|^{2}d\mu(\omega)\Big)^\frac{1}{2}\Big(\int_{\Omega}\|\Lambda_{\omega}h\|^{2}d\mu(\omega)\Big)^\frac{1}{2}\\
&\leq\|\varphi\|_{2}\Big(\int_{\Omega}\|\Lambda_{\omega}h\|^{2}d\mu(\omega)\Big)^\frac{1}{2}\\
&\leq C\|K^*h\|\Big(\int_{\Omega}\|\Lambda_{\omega}h\|^{2}d\mu(\omega)\Big)^\frac{1}{2}.
\end{align*}
Then for each $h\in H$,
$$\frac{1}{C}\|K^*h\|\leq\Big(\int_{\Omega}\|\Lambda_{\omega}h\|^{2}d\mu(\omega)\Big)^\frac{1}{2}.$$
$ (ii)\Rightarrow (iii)$ By Theorem 3.1 in \cite{ckg}, the proof is completed.\\
$ (iiii)\Rightarrow (i)$ 
Assume that there exists a $cg$-Bessel family $\{\Gamma_{\omega}\}_{\omega\in\Omega}$ for $H$ such that
\begin{align}\label{1b}
\langle Kf, h \rangle= \int_{\Omega} \langle  \Gamma_{\omega} f, \Lambda_{\omega} h \rangle d\mu (\omega),\quad  f,h\in H.
\end{align}
So there exists a $C>0$ such that
$$\Big(\int_{\Omega}\|\Gamma_{\omega}f\|^2d\mu (\omega)\Big)^{\frac{1}{2}}\leq C \|f\|, \quad   f\in H.$$
For  $f\in H$, we set $\varphi=\{\Gamma_{\omega}f\}_{\omega\in\Omega}$, then $\varphi\in\big(\oplus_{\omega \in \Omega}H_{\omega},\mu\big)_{L^2}$ and by (\ref{1b}), we obtain
\begin{align*}
\langle Kf, h \rangle= \int_{\Omega} \langle \Lambda_{\omega}^*\varphi(\omega),  h \rangle d\mu (\omega),\quad   f,h\in H.
\end{align*}
Therefore, $\{\Lambda_{\omega}\}_{\omega\in\Omega}$ is an atomic $cg$-system for $K$.
\end{proof}
\begin{theorem}
Let $K_{1}, K_{2}\in B(H)$. If $\{\Lambda_{\omega}\}_{\omega\in\Omega}$ is  an atomic $cg$-system for both $K_{1}$ and $ K_{2}$ and $\alpha, \beta$ are real numbers, then $\{\Lambda_{\omega}\}_{\omega\in\Omega}$ is  an atomic $cg$-system for both operators $\alpha K_{1}+\beta K_{2}$ and $K_{1} K_{2}$.
\end{theorem}
\begin{proof}
By Theorem \ref{3.2}, it is enough to show that $\{\Lambda_{\omega}\}_{\omega\in\Omega}$  is a $c$-($\alpha K_{1}+\beta K_{2}$)-$g$-frame and $c$-$K_{1}K_{2}$-$g$-frame for $H$. Since $\{\Lambda_{\omega}\}_{\omega\in\Omega}$ is an atomic $cg$-system for $K_{1}$ and $K_{2}$, by Theorem \ref{3.2}, $\{\Lambda_{\omega}\}_{\omega\in\Omega}$ is a $c$-$K_{n}$-$g$-frame $(n=1,2)$ for $H$, hence there are positive constants $A_{n}, B_{n}$, $n=1,2$,
such that
\begin{equation}\label{n3.2}
  A_{n}\|K^{*}_{n}f\|^{2}\leq \int_{\Omega} \|\Lambda_{\omega}f\|^{2}\,d\mu(\omega)\leq B_{n}\|f\|^{2}, \quad f\in H.
\end{equation}
For each $\alpha,\beta \in \mathbb{R}$ and $f\in H$, we have
\begin{align*}
\|(\alpha K_{1}^*+\beta K_{2}^*)f\|^{2}\leq 2|\alpha|^{2}\|K_{1}^*f\|^{2}+2|\beta|^{2}\| K_{2}^*f\|^{2},
\end{align*}
therefore
\begin{align*}
\frac{1}{2|\alpha|^{2}|\beta|^{2}}\|(\alpha K_{1}^*+\beta K_{2}^*)f\|^{2} &\leq \frac{1}{|\beta|^{2}}\|K_{1}^*f\|^{2}+\frac{1}{|\alpha|^{2}}\| K_{2}^*f\|^{2}\\
\leq & (\frac{1}{|\beta|^{2}A_1}+\frac{1}{|\alpha|^{2}A_2})\int_\Omega \|\Lambda_\omega f\|^2 d\mu(\omega)\\
=&\frac{|\alpha|^{2} A_2+|\beta|^{2} A_1}{|\alpha|^{2} |\beta|^{2} A_1 A_2}\int_\Omega \|\Lambda_\omega f\|^2 d\mu(\omega).
\end{align*}
Hence
\begin{align*}
\frac{A_1 A_2}{2(|\alpha|^{2} A_2 +|\beta|^{2} A_1)}\|(\alpha K_{1}^*+\beta K_{2}^*)f\|^{2}\leq \int_\Omega \|\Lambda_\omega f\|^2 d\mu(\omega),\quad f\in H.
\end{align*}
Also by inequality $(\ref{n3.2})$, we get
$$\int_{\Omega} \|\Lambda_{\omega}f\|^{2}\,d\mu(\omega)\leq\frac{B_{1}}{2}\|f\|^{2}+\frac{B_{2}}{2}\|f\|^{2}=\frac{B_{1}+B_{2}}{2}\|f\|^{2},\quad f\in H.$$
That is, $\{\Lambda_{\omega}\}_{\omega\in\Omega}$  is a $c$-($\alpha K_{1}+\beta K_{2}$)-$g$-frame for $H$.

\vspace{1.5mm}
Now, for each $f\in H$,
$$\|(K_{1}K_{2})^*f\|^{2}=\|K_{2}^*K_{1}^*f\|^{2}\leq   \|K_{2}^*\|^{2}\|K_{1}^*f\|^{2}.$$
 Since $\{\Lambda_{\omega}\}_{\omega\in\Omega}$ is an atomic $cg$-system for $K_{1}$, for each $f\in H$, we have
$$\frac{A_{1}}{\|K^{*}_{2}\|^{2}} \|(K_{1}K_{2})^*f\|^{2}\leq A_{1}\|K^{*}_{1}f\|^{2}\leq \int_{\Omega} \|\Lambda_{\omega}f\|^{2}\,d\mu(\omega)\leq B_{1}\|f\|^{2}.$$
Therefore $\{\Lambda_{\omega}\}_{\omega\in\Omega}$  is a $c$-$K_{1}K_{2}$-$g$-frame for $H$.
\end{proof}
In the following, we find  some new  atomic $cg$-systems from given ones.
\begin{theorem}
Let $\{\Lambda_{\omega}\}_{\omega\in\Omega}$ and $\{\Gamma_{\omega}\}_{\omega\in\Omega}$ be two atomic $cg$-systems for $K$ and $T_{\Lambda}$ and  $T_{\Gamma}$ be their corresponding  synthesis operators.
Suppose that $T_{\Lambda}T^*_{\Gamma}=0$ and $U, V \in B(H)$ and $U$ is bounded below and $UK^*=K^*U$. Then $\{\Lambda_{\omega}U+\Gamma_{\omega}V\}_{\omega\in\Omega}$  is an atomic $cg$-system for $K$.
\end{theorem}
\begin{proof}
By Theorem \ref{3.2}, we show that $\{\Lambda_{\omega}U+\Gamma_{\omega}V\}_{\omega\in\Omega}$  is a $c$-$K$-$g$-frame for $H$. Since $\{\Lambda_{\omega}\}_{\omega\in\Omega}$ and $\{\Gamma_{\omega}\}_{\omega\in\Omega}$ are   atomic $cg$-systems for $K$, by Theorem \ref{3.2}, $\{\Lambda_{\omega}\}_{\omega\in\Omega}$ and $\{\Gamma_{\omega}\}_{\omega\in\Omega}$ are   $c$-$K$-$g$-frames for $H$  and so there exist $B_{1}\geq A_{1}>0$ and $B_{2}\geq A_{2}>0$ such that for each $f\in H$,
\begin{align*}
A_{1}\|K^{*}f\|^{2}\leq \int_{\Omega} \|\Lambda_{\omega}f\|^{2}\,d\mu(\omega)\leq B_{1}\|f\|^{2},\\
  A_{2}\|K^{*}f\|^{2}\leq \int_{\Omega} \|\Lambda_{\omega}f\|^{2}\,d\mu(\omega)\leq B_{2}\|f\|^{2}.
\end{align*}
Since $T_{\Lambda}T^*_{\Gamma}=0$, for each $f\in H$, we have
$$\int_{\Omega}\langle\Lambda_{\omega}^{*}\Gamma_{\omega}f, f \rangle d\mu(\omega)=0.$$
Therefore, for each $f\in H$, we can write
\begin{align*}
\int_{\Omega}\Vert(\Lambda_{\omega}U+\Gamma_{\omega}V)\Vert^2 d\,\mu(\omega)
&=\int_{\Omega}\Vert\Lambda_{\omega}Uf\Vert^2 d\,\mu(\omega)+\int_{\Omega}\Vert\Gamma_{\omega}Vf\Vert^2 d\,\mu(\omega)\\
&\leq B_{1}\|Uf\|^{2}+B_{2}\|Vf\|^{2}\\
&\leq(B_{1}\|U\|^{2}+B_{2}\|V\|^{2})\|f\|^{2}.
\end{align*}
That is, $\{\Lambda_{\omega}U+\Gamma_{\omega}V\}_{\omega\in\Omega}$  is a $cg$-Bessel family for $H$. Now, we show that $\{\Lambda_{\omega}U+\Gamma_{\omega}V\}_{\omega\in\Omega}$ has the lower $c$-$K$-$g$-frame condition.
Since $U$ is bounded below, so there exists a $C>0$ such that for each $f\in H$,  $\|Uf\|^{2}\geq C\|f\|^{2}$. By  the assumption, $UK^*=K^*U$, so we have
\begin{align*}
\int_{\Omega}\Vert(\Lambda_{\omega}U+\Gamma_{\omega}V)f\Vert^2 d\,\mu(\omega)
&=\int_{\Omega}\Vert\Lambda_{\omega}Uf\Vert^2 d\,\mu(\omega)+\int_{\Omega}\Vert\Gamma_{\omega}Vf\Vert^2 d\,\mu(\omega)\\
&\geq\int_{\Omega}\Vert\Lambda_{\omega}Uf\Vert^2 d\,\mu(\omega)\\
&\geq A_{1}\|K^*Uf\|^{2}=A_{1}\|UK^*f\|^{2}\\
&\geq CA_{1}\|K^*f\|^{2}.
\end{align*}
So  $\{\Lambda_{\omega}U+\Gamma_{\omega}V\}_{\omega\in\Omega}$ is a $c$-$K$-$g$-frame for $H$ and by Theorem \ref{3.2}, it is an atomic $cg$-system for $K$.
\end{proof}
\begin{corollary}
 Let  $K\in B(H)$ and $\{\Lambda_{\omega}\}_{\omega\in\Omega}$ be an atomic $cg$-system for $K$. If $U\in B(H)$ is bounded below operator and $UK^*=K^*U$, then $\{\Lambda_{\omega}U\}_{\omega\in\Omega}$  is an atomic $cg$-system for $K$.
\end{corollary}

If $U=V=I$, then we have the following result:

\begin{corollary}
 Let  $\{\Lambda_{\omega}\}_{\omega\in\Omega}$  and $\{\Lambda_{\omega}\}_{\omega\in\Omega}$ be two Parseval $c$-$K$-$g$-frame for $H$, with synthesis operators $T_{\Lambda}$ and  $T_{\Gamma}$, respectively. If $T_{\Lambda}T^*_{\Gamma}=0$, then $\{\Lambda_{\omega}+\Gamma_{\omega}\}_{\omega\in\Omega}$  is a $2$-tight $c$-$K$-$g$-frame for $H$.
\end{corollary}
\begin{theorem}
Let $(\Omega, \mu)$ be a measure space, where $\mu$ is $\sigma$-finite. Suppose that $\{\Lambda_{\omega}\}_{\omega\in\Omega}$ and $\{\Gamma_{\omega}\}_{\omega\in\Omega}$ are  atomic $cg$-systems for $K$  and $T_{\Lambda}$ and  $T_{\Gamma}$ are the  synthesis operators of $\{\Lambda_{\omega}\}_{\omega\in\Omega}$ and $\{\Gamma_{\omega}\}_{\omega\in\Omega}$, respectively. If $T_{\Lambda}T^*_{\Gamma}=0$ and $U_1, U_2\in B(H)$ are such that  $R(T_{\Lambda})\subseteq R(U^*_{1}T_{\Lambda})$, $R(T_{\Gamma})\subseteq R(U^*_{2}T_{\Gamma})$, then $\{\Lambda_{\omega}U_{1}+\Gamma_{\omega}U_{2}\}_{\omega\in\Omega}$  is an atomic $cg$-system for $K$.
\end{theorem}
\begin{proof}
Since $T_{\Lambda}T^*_{\Gamma}=0$, for each $f\in H$, we have
\begin{align}\label{as1}
\int_{\Omega}\Vert(\Lambda_{\omega}U_{1}+\Gamma_{\omega}U_{2})f\Vert^2 d\,\mu(\omega)&=\int_{\Omega}\Vert\Lambda_{\omega}U_{1}f\Vert^2 d\,\mu(\omega)+\int_{\Omega}\Vert\Gamma_{\omega}U_{2}f\Vert^2 d\,\mu(\omega)  \notag  \\
&=\|T^*_{\Lambda}U_{1}f\|^{2}_2+\|T^*_{\Gamma}U_{2}f\|^{2}_2  \notag \\
&=\|(U^*_{1}T_{\Lambda})^*f\|^{2}_2+\|(U^*_{2}T_{\Gamma})^*f\|^{2}_2.
\end{align}
Since $\{\Lambda_{\omega}\}_{\omega\in\Omega}$ and $\{\Gamma_{\omega}\}_{\omega\in\Omega}$ are atomic $cg$-systems for $K$, by Theorem \ref{3.2}, they are $c$-$K$-$g$-frames for $H$. So by Theorem 2.5 in \cite{ckg}, we have
$$R(K)\subseteq R(T_{\Lambda})\subseteq R(U^*_{1}T_{\Lambda})$$
and
$$R(K)\subseteq R(T_{\Gamma})\subseteq R(U^*_{2}T_{\Gamma}).$$
Thus by Lemma \ref{1.1},  there exist $\lambda_{1}, \lambda_{2}>0$ such that
$$KK^*\leq\lambda_{1}(U^*_{1}T_{\Lambda})(U^*_{1}T_{\Lambda})^*$$
and
$$KK^*\leq\lambda_{2}(U^*_{2}T_{\Gamma})(U^*_{2}T_{\Gamma})^*.$$
By (\ref{as1}), for each $f\in H$, we have
\begin{align*}
\int_{\Omega}\Vert(\Lambda_{\omega}U_{1}+\Gamma_{\omega}U_{2})\Vert^2 d\,\mu(\omega)&=\|(U^*_{1}T_{\Lambda})^*f\|^{2}_2+\|(U^*_{2}T_{\Gamma})^*f\|^{2}_2\\
&\geq(\frac{1}{\lambda_{1}}+\frac{1}{\lambda_{2}})\|K^*f\|^{2}.
\end{align*}
Hence $\{\Lambda_{\omega}U_{1}+\Gamma_{\omega}U_{2}\}_{\omega\in\Omega}$  is a $c$-$K$-$g$-frame for $H$ and by Theorem \ref{3.2}, the conclusion holds.
\end{proof}

\begin{theorem}
 Let  $K\in B(H)$ and $K$ be with closed range. Suppose  that $\{\Lambda_{\omega}\}_{\omega\in\Omega}$ is an atomic $cg$-system for $K$ and $S_{\Lambda}$ is the frame operator of $\{\Lambda_{\omega}\}_{\omega\in\Omega}$. If $U\in B(H)$ is a positive operator such that $US_{\Lambda}=S_{\Lambda}U$, then $\{\Lambda_{\omega}+\Lambda_{\omega}U\}_{\omega\in\Omega}$ is an atomic $cg$-system for $K$. Moreover, for each $n\in \mathbb{N}$, $\{\Lambda_{\omega}+\Lambda_{\omega}U^{n}\}_{\omega\in\Omega}$ is an atomic $cg$-system for $K$.
\end{theorem}
\begin{proof}
By Theorem \ref{3.2},   $\{\Lambda_{\omega}\}_{\omega\in\Omega}$ is a $c$-$K$-$g$-frame  for $H$. Then by Lemma \ref{x}, there exists $A>0$ such that $S_{\Lambda}\geq AKK^{*}$. The frame  operator of $\{\Lambda_{\omega}+\Lambda_{\omega}U\}_{\omega\in\Omega}$ is given by
$$(I+U)^*S_{\Lambda}(I+U).$$
By assumption, $K$ has closed range, thus $S_{\Lambda}$ is positive and
$$(I+U)^*S_{\Lambda}(I+U)\geq AKK^*.$$
By Lemma \ref{x}, we conclude  that $\{\Lambda_{\omega}+\Lambda_{\omega}U\}_{\omega\in\Omega}$ is a $c$-$K$-$g$-frame for $H$ and by Theorem \ref{3.2}, it is an atomic $cg$-system for $K$. For  any $n\in \mathbb{N}$, the frame operator of $\{\Lambda_{\omega}+\Lambda_{\omega}U^{n}\}_{\omega\in\Omega}$
 is $(I+U^{n})^*S_{\Lambda}(I+U^{n})$ and similarly $\{\Lambda_{\omega}+\Lambda_{\omega}U^{n}\}_{\omega\in\Omega}$ is an atomic $cg$-system for $K$.
\end{proof}
\begin{theorem}
 Let  $K\in B(H)$ and $K$ be with closed range, $\{\Lambda_{\omega}\}_{\omega\in\Omega}$ be  an atomic $cg$-system for $K$. Assume that $\{\Gamma_{\omega}\}_{\omega\in\Omega}$ is a $cg$-Bessel family with respect to $\{H_{\omega}\}_{\omega\in\Omega}$ such that
\begin{align*}\label{m}
\langle Kf, h \rangle =\int_{\Omega}\langle\Lambda_{\omega}^{*}\Gamma_{\omega}f, h \rangle d\mu(\omega), \quad f, h\in H.
\end{align*}
Then the family $\{\Theta_{\omega}\}_{\omega\in\Omega}=\{\Gamma_{\omega}U\}_{\omega\in\Omega}$, where $U=K^{\dagger}|_{R(K)}$, is a $cg$-frame for $R(K)$  such that
$$\langle f, h \rangle =\int_{\Omega}\langle\Lambda_{\omega}^{*}\Theta_{\omega}f, h \rangle d\mu(\omega)=\int_{\Omega}\langle\Theta_{\omega}^{*}\Lambda_{\omega}f, h \rangle d\mu(\omega), \quad f\in R(K), h\in H.$$
\end{theorem}
\begin{proof}
Since $\{\Lambda_{\omega}\}_{\omega\in\Omega}$ is an atomic $cg$-system for $K$, by  Theorem \ref{3.2} (iii), there exists a $cg$-Bessel family $\{\Gamma_{\omega}\}_{\omega\in\Omega}$ with bound $B$  such that
\begin{equation}\label{au1}
\langle Kf, h \rangle =\int_{\Omega}\langle\Lambda_{\omega}^{*}\Gamma_{\omega}f, h \rangle d\mu(\omega), \quad f, h\in H.
\end{equation}
By assumption, $R(K)$ is closed, so there exists Pseudo-inverse $K^\dagger$ of $K$, such that
$$f=KK^\dagger f, \quad f\in R(K).$$
From (\ref{au1}), we have
$$\langle f, h \rangle =\langle KK^\dagger f, h \rangle =\int_{\Omega}\langle\Lambda_{\omega}^{*}\Gamma_{\omega}K^\dagger f, h \rangle d\mu(\omega)~,~f\in R(K), h\in H$$
Now, let $\Theta_{\omega}=\Gamma_{\omega}U$ where $U=K^{\dagger}|_{R(K)}$, so  $\Theta_{\omega}:R(K)\longrightarrow H_{\omega}.$
For each $f\in R(K), ~K^\dagger f \in H$  and
$$\int_{\Omega}\|\Theta_{\omega}f\|^{2}d\mu(\omega)=\int_{\Omega}\|\Gamma_{\omega}K^\dagger f\|^{2}d\mu(\omega)\leq B\|K^\dagger f\|^{2}\leq B\|K^\dagger \|^{2}\|f\|^{2}.$$
That is, $\{\Theta_{\omega}\}_{\omega\in\Omega}$ is a $cg$-Bessel family for $R(K)$ with respect to $\{H_{\omega}\}_{\omega\in\Omega}$. Now, we show that $\{\Theta_{\omega}\}_{\omega\in\Omega}$ has the lower frame condition. For each $ f \in R(K)$,
\begin{align*}
\|f\|^{2}=|\langle f, f \rangle |&=|\langle KK^\dagger f, f \rangle |=\Big|\int_{\Omega}\langle\Lambda_{\omega}^{*}\Gamma_{\omega}K^\dagger f, f \rangle d\mu(\omega)\Big|\\
&\leq \Big|\int_{\Omega}\langle \Gamma_{\omega}Uf, \Lambda_{\omega}f \rangle d\mu(\omega)\Big|\\
&\leq \Big(\int_{\Omega}\|\Gamma_{\omega}Uf\|^{2}d\mu(\omega)\Big)^{\frac{1}{2}}\Big(\int_{\Omega}\|\Lambda_{\omega}f\|^{2}d\mu(\omega)\Big)^{\frac{1}{2}}\\
&\leq C\|f\|\Big(\int_{\Omega}\|\Gamma_{\omega}Uf\|^{2}d\mu(\omega)\Big)^{\frac{1}{2}},
\end{align*}
where $C$ is the upper frame bound of $c$-$K$-$g$-frame  $\{\Lambda_{\omega}\}_{\omega\in\Omega}$.\\
Therefore
$$\frac{1}{C}\|f\|\leq \Big(\int_{\Omega}\|\Gamma_{\omega}Uf\|^{2}d\mu(\omega)\Big)^{\frac{1}{2}}, \quad f\in R(K).$$
Thus $\{\Theta_{\omega}\}_{\omega\in\Omega}=\{\Gamma_{\omega}U\}_{\omega\in\Omega}$  is a $cg$-frame for $R(K)$.\\
The rest of proof is straightforward.
\end{proof}
\bibliographystyle{amsplain}

\begin{thebibliography}{10}

\bibitem{1} {\it M. R. Abdollahpour and M. H. Faroughi}, Continuous $G$-frames in Hilbert spaces,
Southeast Asian Bull. Math., {\bf 32} (2008), 1-19.

\bibitem{nob} {\it E. N. Alamdari, M. Azadi and H. Doostie}, Excess of continuous $K$-$g$-frames and some other properties,
Submitted.

\bibitem{2} {\it S. T. Ali, J. P. Antoine and J. P. Gazeau}, Continuous frames in Hilbert spaces,
Annals of Phy., { \bf 222 }(1993), 1-37.

\bibitem{ckg} { \it E. Alizadeh, A. Rahimi, E. Osgooei and M. Rahmani}, Continuous $K$-$g$-frames in Hilbert spaces,
Bull. Iran. Math. Soc., {\bf 45} (4) (2019), 1091-1104

\bibitem{sckg} { \it E. Alizadeh, A. Rahimi, E. Osgooei and M. Rahmani},  Some properties of Continuous $K$-$G$-frames in Hilbert spaces,
U. P. B. Sci. Bull, Series A.,  {\bf 81} (3) (2019), 43-52.

\bibitem{ar} {\it F. Arabyani Neyshaburi, G. Mohajeri Minaei and E. Anjidani}, On some equalities and inequalities for $K$-frames,
 Indian. J. Pure. Appl. Math., (2018).

\bibitem{dag}{\it R. G. Douglas}, On majorization, factorization and range inclusion of operators on Hilbert space,
Pro. Amer. Math. Sco., {\bf 17} (2) (1966), 413-415.

\bibitem{Duff}{\it R. J. Duffin and A. C. Schaeffer}, A class of nonharmonic Fourier series,
Trans. Amer. Math. Soc., {\bf 72} (1952), 341-366.

\bibitem{k-frame}{\it L. Gavruta}, Frames for operators,
Appl. Comput. Harmon. Anal., {\bf 32} (2012), 139-144.

\bibitem{10}{\it D. Hua and Y. Hung}, $K$-$g$-frames and stability of $K$-$g$-frames in Hilbert spaces,
J. Korean Math. Soc., {\bf 53} (6) (2016), 1331-1345.

\bibitem{Li}{\it D. Li, J. Leng and T. Huang}, Generalized frames for operators associated with atomic systems,
Banach J. Math. Anal., {\bf 12} (1) (2018), 206-221.

\bibitem{14}{\it A. Rahimi, A. Najati and Y. N. Dehghan}, Continuous frames in Hilbert spaces,
Method. Funct. Anal. Topology.,  {\bf 12} (2) (2006), 170-182.

\bibitem{Rah}{\it M. Rahmani}, Sums of $c$-frames, $c$-Riesz Bases and orthonormal mappings,
U. P. B. Sci. Bull, Series A.,  {\bf 77} (3) (2015), 3-14.

\end{thebibliography}

\vspace{7mm}
\end{document}